\input amstex
\documentstyle{amsppt}

\NoRunningHeads \magnification=1200 \TagsOnRight \NoBlackBoxes
\hsize=5.3in \vsize=7.2in \hoffset= 0in \baselineskip=12pt
\topmatter
\title  Local structure of closed  symmetric 2-differentials
\endtitle
\rightheadtext{Closed symmetric differentials on surfaces}
\author Fedor Bogomolov*
 Bruno De Oliveira**
\endauthor
\thanks
* Partially supported by  AG Laboratory GU-HSE grant RF government ag. 1111.G34.31.0023.
** The second author thanks the financial support  by the Courant Institute of Mathematical Sciences during 2011/12.
\endthanks
\affil
* Courant Institute for Mathematical Sciences
**University of Miami
\endaffil
\address
Fedor Bogomolov Courant Institute for Mathematical Sciences,
New York University\\
Bruno De Oliveira University of Miami
\endaddress
\email bogomolo{\@}CIMS.NYU.EDU bdeolive@math.miami.edu
\endemail
\endtopmatter
\document

\head {0. Introduction}\endhead

In the authors's previous work on symmetric differentials and their connection to  the topological properties
of the ambient manifold, a class of symmetric differentials was introduced: closed symmetric differentials
([BoDeO11] and [BoDeO13]). Closed symmetric differentials are characterized by the possibility to locally
decompose the differential as a product of closed holomorphic 1-differentials  in a neighborhood of a  point
of the manifold. The property of being closed is conjecturally described by a  non-linear differential operator
(in the case of dimension 2 and degree 2 this differential operator comes from the Gaussian curvature, see section 2.1).

\

In this article we give a description of the local structure of closed symmetric 2-differentials on complex surfaces,
with an emphasis towards the local decompositions as products of 1-differentials. Recall  that there is a general
obstruction for a symmetric 2-differential  to have a decomposition  as a product of 1-differentials around a point
$x$ in the complex surface $X$, it might be impossible to order the two foliations defined by $w$ near $x$ (we then
say that $w$ is not locally split at $x$). This obstruction can be removed via a ramified covering of $X$, hence the
results will be given for symmetric differentials that are locally split.

\

We show that a closed symmetric 2-differential $w$ of rank 2 (i.e. defines two distinct foliations at the general point)
has a  subvariety $B_w\subset X$ outside of which  $w$ is locally the  product  of closed holomorphic 1-differentials.
The main result, theorem 2.6, gives a complete description of a (locally split) closed symmetric 2-differential in a
neighborhood of  a general point of $B_w$. A consequence of the main result is that the  differential $w$ still has
a local decomposition  into a product of closed 1-differentials (in a generalized sense) at the points of $B_w$.
The closed 1-differentials involved in the local decompositions might have to be multi-valued and acquire singularities
along $B_w$. Note that if we were considering local decompositions of  a locally split holomorphic symmetric 2-differential
into a product of 1-differentials (not necessarily closed), then  the 1-differentials involved can be chosen to be
holomorphic, i.e. no singularities need to occur. On the other hand, it is also true that by multiplying one of the
1-differentials by an arbitrary function and the other by its inverse, that arbitrary singularities can occur in the
decomposition. An important feature of decompositions of symmetric 2-differentials of rank 2 as products of closed
1-differentials is that they are unique up to multiplicative constants, hence there is no  ambiguity on the singularities that occur.

\

The singularities that occur in the decomposition of a closed
holomorphic symmetric 2-differential $w$ when we require that the
1-differentials  are closed can be essential singularities along
the locus $B_w$. A key feature of  theorem 2.6 giving the local
structure of $w$ around points in $B_w$ is that we  have a control
on these  essential singularities, they come from exponentials of
meromorphic functions acquiring poles of a bounded order along
$B_w$. Before describing the nature of the bound, we need to
describe our result characterizing the locus $B_w$. In the case
$w$ is locally split (always the case after a ramified cover), we
show that any irreducible component of $B_w$ must be
simultaneously a leaf of both foliations defined by $w$. The bound
on the order of the poles along an irreducible component of $B_w$
is the order of contact of both foliations along that irreducible
component.

\

This article addresses the case of closed symmetric 2-differentials, we expect that a straightforward generalization of our methods will
provide similar results on the local structure of closed symmetric differentials of arbitrary degree and give control of the singularities that occur on  the decompositions as product of closed 1-differentials.

  \

  \

\head {1. General set up}\endhead

\

\

A symmetric differential $w\in H^0(X,S^m\Omega_X^1)$ on a complex manifold $X$ defines at each point where $w(x)\neq 0$ a cone in tangent space $T_xX$ with vertex the origin and defined at infinity ($\Bbb P^{n-1}$) by a variety of degree m. If $X$ is a complex surface then one gets a  distribution of $d$ $(d\le m)$ lines, which will be integrable, defining a non-singular d-web at the general point. In higher dimensions the cones will not be necessarily union of hyperplanes and even if they are hyperplanes their distributions need not be integrable. Here, we should note that the class of symmetric differentials that is studied in this work, closed symmetric differential (see below),  will in all dimensions be connected to webs on the manifold.

\proclaim {Definition 1.1} A symmetric differential $w\in H^0(X,S^m\Omega^1_X)$ is split if it has a decomposition:

$$w=\psi_1...\psi_m$$

\noindent where the $\psi_i$ are meromorphic 1-forms or equivalently if
$w=\mu_1...\mu_m$, with $\mu_i\in H^0(X,\Omega^1_X\otimes L_i)$,  where $L_i$ are line bundles on $X$.
\endproclaim

Geometrically being split means that the symmetric differential defines hyperplane distributions  and moreover they can be numbered consistently globally.

\proclaim {Definition 1.2} A symmetric differential $w$ on $X$ has rank r if at a general point $x\in X$ $w(x)$ defines $r$ distinct hyperplanes in $T_xX$.
\endproclaim

\proclaim {Definition 1.3} A symmetric differential $w$ on $X$ is said to have a holomorphic closed decomposition if:

$$w=\mu_1...\mu_m \text { },\text { }\text { }\text { }\text { }\text { }\text { }\text { }\text { }\text { }\mu_i \text { closed holomorphic 1-forms} \tag 1.1$$

\noindent and a holomorphic closed decomposition at $x$ if $x$ has an analytic  neighborhood where (1.1) holds.
\endproclaim

\proclaim {Definition 1.4} A symmetric differential $w \in H^0(X,S^m\Omega^1_X)$ is said to be:
\smallskip
1) closed, if $w$ has  an holomorphic closed decomposition at a general point $x\in X$.
\smallskip
2) of the 1st kind, if $w$ has holomorphic closed decompositions at all $x\in X$.\endproclaim

\

\noindent Remarks: 1) The class of closed symmetric differentials of the 1st kind plays a special role in the motivation for considering closed symmetric differentials as a class of symmetric differentials having a stronger connection to the topology of the ambient manifold. We expand on this point below.
\smallskip
2) Our definitions of closed and 1st kind coincide with the usual definitions when m=1, i.e. holomorphic 1-forms. Our definition of closed asks for a holomorphic 1-form to be locally exact somewhere which by the identity principle implies it is locally exact everywhere and hence closed in the usual sense. Hence, for $m=1$ our notions of closed and 1st kind  coincide.
\smallskip 3) If the degree $m>1$, then closed no longer implies of the 1st kind. This has far reaching geometric consequences and the main results of this work concern  the locus where this failure comes from and the structure of the closed symmetric differentials near this locus.

\proclaim {Definition 1.5} The locus of $X$ where a closed symmetric differential $w$ fails to be of the 1st kind at, $B_w=\{x\in X|\text { }w \text { has no holomorphic closed decomposition at x}\}$, will be called the breakdown locus of $w$.
\endproclaim

A key feature of holomorphic closed decompositions is that they have rigidity properties. The level of rigidity has to do with a familiar notion in the theory of webs, the abelian rank of a web.

\proclaim {Definition 1.6} Given the germ $w_x\in S^m\Omega^1_{X,x}$ with the holomorphic closed decomposition

$$w_x=\mu_1...\mu_m \tag 1.2$$

\noindent where $\mu_i\in \Omega^1_{X,cl,x}$ ($\Omega^1_{X,cl}$ is the sheaf of closed holomorphic 1-forms on $X$), we  call  an m-tuple $(f_1,...,f_m)\in \Cal M_x^m$ satisfying:

$$\sum_{i=1}^mf_i\mu_i=0 \text { }\text { }\text { }\text { } with \text { }\text { }  df_i\wedge \mu_i\equiv 0.$$

\noindent an abelian relation of the decomposition (1.2). The abelian rank of the decomposition (1.2) is the dimension of the $\Bbb C$-vector space consisting of all abelian relations of (1.2). The abelian rank of a closed symmetric differential $w\in H^0(X,S^m\Omega^1_X)$ is the abelian rank of any holomorphic closed  decomposition at the general point of $x$.
\endproclaim

\

\noindent Remarks: 1)  The definition of abelian rank of a closed symmetric differential $w$ is well defined, since there is an analytic subvariety of $R\subset X$  such that all holomorphic closed decompositions of $w_x$, $\forall x\in X\setminus R$, have the same abelian rank.
\smallskip
2) It is a classical result of web theory  that the abelian rank of  a decomposition (1.2) is finite if rank(w)=m, with  upper bounds depending on the dimension and degree (for dimension 2 this is a result of G.Bol and also W.Blaschke see for example [ChGr78], [He01] and [He04] for information on webs).
\smallskip

3)  A general germ of a closed symmetric differential has trivial abelian rank.

\

We concentrate our attention to the case of trivial abelian rank which is the generic case and  holds trivially for all closed symmetric 2-differentials (and rank 2).

\proclaim {Proposition 1.1} Let $X$ be a connected manifold and $w\in H^0(X,S^m\Omega^1_X)$ be a closed symmetric differential with an holomorphic  closed decomposition $w=\mu_1...\mu_m$. If the abelian rank of $w$ is trivial, then all holomorphic closed decomposition $w=\eta_1...\eta_m$ of $w$ on $X$ have the closed 1-forms $\eta_i=c_i\mu_i$ with $c_i\in \Bbb C^*$ and $\prod_{i=1}^mc_i=1$.
\endproclaim

\demo {Proof} Suppose $w=\eta_1...\eta_m$ is an holomorphic  closed decomposition of $w$ and assume the $\eta_i$ are ordered such that $\eta_i\wedge\mu_i=0$. The condition $\eta_i\wedge\mu_i=0$ in conjunction with $\eta_i$ and $\mu_i$ being closed implies that $\eta_i=f_i\mu_i$ with $f_i\in \Cal M(X)$ and $df_i\wedge \mu_i=0$. Moreover $\mu_1...\mu_m=\eta_1...\eta_m$ gives:

$$\prod_{i=1}^mf_i=1\tag 1.3$$

\noindent Pick a simply connected open set $U\subset X$ where $f_i|_U\in \Cal O^*(U)$. Taking the logarithm and differentiating (1.3) restricted to $U$ and using the identity principle we obtain:

$$\sum_{i=1}^m\frac{df_i}{f_i}=0 \tag 1.4$$

\noindent but $df_i\wedge \mu_i=0$, hence $df_i=g_i\mu_i$ with $g_i\in \Cal M(X)$ and $dg_i\wedge \mu_i=0$. It follows that (1.4) gives rise to the abelian relation at a general point $x\in X$:

$$\sum_{i=1}^m(\frac{g_i}{f_i})_x(\mu_i)_x=0$$

\noindent The abelian rank of $w$ being trivial implies that the $(g_i)_x=0$ and hence $df_i=0$ on $X$, i.e. $f_i=c_i\in \Bbb C^*$.$\blacklozenge$
\enddemo

\

A consequence of this proposition, see below, is the decomposition of symmetric differentials of the 1st kind with abelian rank 0 into a product of twisted closed holomorphic 1-forms $\phi_i\in H^0(X,\Omega^1_{X,cl}\otimes \Bbb C_{\rho_i})$. In [BoDeO13] we use such decompositions to characterize the origins and the geometric implications of symmetric 2-differentials of the 1st kind.

\proclaim {Corollary 1.2} Let  $X$ be a complex manifold and $w\in H^0(X,S^m\Omega^1_X)$  be of the 1st kind with trivial abelian rank. Then  there is a finite unramified cover $f:X'\to X$ (unnecessary if $w$ is split) for which $f^*w$ has a decomposition:

$$f^*w=\phi_1...\phi_m$$

\noindent where $\phi_i\in H^0(X',\Omega^1_{X',cl}\otimes \Bbb C_{\rho_i})$, where the $\Bbb C_{\rho_i}$ are local systems of rank 1 on $X'$ such that $\Bbb C_{\rho_1}\otimes ... \otimes \Bbb C_{\rho_m}\simeq \Bbb C$.
\endproclaim

\demo {Proof} The differential $w$ being of first kind implies that  locally $w$ is split, but $w$ might fail to be globally split. This failure is measured by the monodromy coming from the local ordering of the foliations, i.e. we obtain a representation $\sigma:\pi_1(X,x) \to S_m$. Associated to this representation we get an unramified cover $f:X' \to X$ with degree a factor of $m!$ such that $f^*w$ is split.

\

From now on we assume that $w$ is split on $X$. The differential $w$ being of the 1st kind gives that there is a Leray covering $\Cal U=\{U_i\}$ of $X$ where

$$w|_{U_i}=\phi_{1i}...\phi_{mi}$$

\noindent with $\phi_{ki}\in H^0(U_i,\Omega^1_{X,cl})$. Since the differential $w$ is split, we can order the $\{\phi_{ki}\}$ such that $\phi_{ki}\wedge \phi_{kj}=0$ on the intersections $U_i\cap U_j$. Proposition 1.1 implies that on $U_i\cap U_j$

$$\phi_{ki}=c_{k,ij}\phi_{kj} \tag 1.5$$

\noindent with $c_{k,ij}\in \Bbb C^*$ and $\prod_{k=1}^mc_{k,ij}=1$. The $m$ collections $\{c_{k,ij}\}$ for $k=1,...,m$ are elements in $Z^1(X,\Bbb C^*)$ and give rise to the rank 1 local systems which we denote by $\Bbb C_{\rho_k}$ and satisfy $\Bbb C_{\rho_1}\otimes ... \otimes \Bbb C_{\rho_m}\simeq \Bbb C$ (we remark that the isomorphism classes of these local systems are completely determined by $w$). It follows from (1.5) that each collection for a fixed $k$, $\{\phi_{ik}\}$, gives a section $\phi_k \in H^0(X,\Omega^1_{X,cl}\otimes \Bbb C_{\rho_k})$ and the result holds.$\blacklozenge$
\enddemo

\

\noindent  The presence of twisted closed holomorphic 1-forms $\phi_i\in H^0(X,\Omega^1_{X,cl}\otimes \Bbb C_{\rho_i})$ has implications on both the topology and geometry of the manifold $X$. On the topological side one observes that the cohomology exact sequence  associated to the short exact sequence $0\to \Bbb C_{\rho} \to \Cal O\otimes \Bbb C_{\rho} \to \Omega^1_{X,cl}\otimes \Bbb C_{\rho} \to 0$ implies that $H^1(X,\Bbb C_{\rho})\ge h^0(X,\Omega^1_{X,cl}\otimes \Bbb C_{\rho})$ (hence in particular $\pi_1(X)$ must infinite). On the geometric side, if $X$ is compact K\"{a}hler the presence of non-torsion, i.e. $L_\rho=\Cal O\otimes \Bbb C_\rho$ non-torsion, twisted closed holomorphic 1-forms  implies that $X$ is fibered   over curves of genus $g\ge 1$ as follows from the work of Beauville, Lazarsfeld-Green and Simpson (see [GrLa87],[Be92],[Si93] and [Ar92]).

\

\

\

\head {2. Local structure of closed 2-differentials on surfaces}\endhead

\

\

\noindent  {\bf 2.1 Differential operator for closed symmetric 2-differentials}

\

A symmetric differential of degree 2 on a complex surface can be viewed as a generalized  complex
counterpart of a Riemannian metric on a real surface. We are going to use this fact to motivate the
differential operator characterizing closed symmetric 2-differentials.

\

Let $w\in H^0(X,S^2\Omega^1_X)$ on a complex surface $X$. There is an open cover of $X$,
$\Cal U=\{U_i\}$ by local holomorphic charts, where

$$w|_{U_i}=a_i(z)(dz_{1i})^2+b_i(z)dz_{1i}dz_{2i}+c_i(z)(dz_{2i})^2$$

\noindent On each of these open sets we get the holomorphic
functions $det w|_{U_i}=a_i(z)c_i(z)-b_i(z)^2/4$ which together
form an element of $H^0(X,\Cal O(2K_X))$, called the discriminant
of $w$.

\proclaim {Definition 2.1} The discriminant divisor $\text
{Disc}_w$ of $w\in H^0(X,S^2\Omega^1_X)$ is the divisor of zeros
of the section $\{det w|_{U_i}\}_{i\in I}$ of $\Cal O(2K_X)$. The
core discriminant divisor of $w$ is $ {\text {Disc}}^0_w=\text
{Disc}_w-2(w)_0$.
\endproclaim

\noindent Geometrically, the support of $\text {Disc}_w$
corresponds to the set of points where $w(x)$ either vanishes or
defines one single line in $T_xX$. To better understand the
support of $\text {Disc}^0_w$ we give a new characterization of
the divisor $\text {Disc}^0_w$. Let $\Cal U=\{U_i\}_{i\in I}$ be
an open covering of $X$ such that $w|_{U_i}=h_i\hat {w_i}$ with
$h_i\in \Cal O(U_i)$ and $\hat {w_i}$ vanishes only in codimension
2. The divisor $\text {Disc}^0_w$ can described via the local
information $\text {Disc}^0_w\cap U_i=\text {Disc}_{\hat {w_i}}$.
By definition $\text {Supp} (\text {Disc}^0_w)\subset \text {Supp}
(\text {Disc}_w)$, moreover a point $x\in \text {Supp} (\text
{Disc}_w)\setminus \text {Supp} (\text {Disc}^0_w)$ is such that
$w(x)=0$ but $\hat {w_i}(x)$ defines two distinct lines in $T_xX$.
A more detailed characterization of the irreducible components of
the divisors $\text {Disc}_w$ and $\text {Disc}^0_w$ is given in
section 2.2.

\

Let $x\in X\setminus \text {Disc}^0_w$, then the germ $w_x$ of $w$
at $x$ splits, $w_x=\mu_1\mu_2$ with $\mu_i\in \Omega^1_{X,x}$.
Say $x\in U_i$ with $U_i$ as above, $w|_{U_i}=h_i\hat {w_i}$, then
$x\in X\setminus \text {Disc}^0_w$ implies that the discriminant
of $\hat {w_i}$ does not vanish at $x$ and therefore $\hat {w_i}$
(and hence $w$)  splits at $x$, $\hat {w_i}_x=\hat \mu_1\hat
\mu_2$. Moreover, since the discriminant of $\hat {w_i}$ is
nonzero at $x$, then $\hat \mu_1\wedge \hat \mu_2(x)\neq 0$, which
implies that $x$ has a neighborhood $U_x$ with holomorphic a chart
$(z_1,z_2)$ such that:

$$w|_{U_x}=g(z)dz_1dz_2 \tag {2.1}$$

\noindent with $g\in \Cal O(U_x)$. The condition that $w|_{U_x}$
has a closed holomorphic decomposition, $w|_{U_x}=g(z)
dz_1dz_2=\mu_1\mu_2$ with $\mu_i$ closed holomorphic 1-forms, is
equivalent to

$$g(z)=f_1(z_1)f_2(z_2)  \tag {2.2}$$

\noindent ($\mu_i\wedge dz_i\equiv 0$ implies that
$\mu_i=f_i(z_i)dz_i$). The condition (2.2) can be characterized
via the nonlinear differential   equation $\frac {\partial ^2 ln
\text { }g(z)}{\partial z_1\partial z_2}=0$.

\smallskip

It follows from the Brioschi's formula for the Gaussian curvaure
in terms of the 1st fundamental form [Sp99], that the differential
operator giving the Gaussian curvature for a metric in the form
$ds^2=f(x)dx_1dx_2$ is $K(ds^2)=-\frac {2}{f} \frac {\partial ^2
ln \text { }f(z)}{\partial x_1\partial x_2}$. Hence the symmetric
differential in (2.1) is closed if and only if $K_\Bbb
C(w|_{U_x})=0$, where $K_\Bbb C$ is the operator obtained from $K$
by replacing $x_1,x_2$ by $z_1,z_2$.

\

A symmetric 2-differential can not be put locally in the form
(2.1) everywhere, but this is not a problem since the differential
operator $K$ for the Gaussian curvature works for metrics whose
1st fundamental form is arbitrary (works formally if the 1st
fundamental form is degenerate, i.e. with discriminant zero at
some points). Hence $K_\Bbb C$ works for any symmetric
2-differential (for a general symmetric 2-differential $w$,
$K_\Bbb C(w)$ will be a meromorphic function with poles along the
discriminant locus). If $w$ is a symmetric 2-differential
satisfying $K_{\Bbb C}(w)=0$, then $K_{\Bbb C}(w|_{U_x})=0$ with
$U_x\subset X\setminus \text {Disc}^0_w$ as in 2.1 and hence it by
our previous paragraph $w$ is closed. Hence we obtain:

\proclaim {Proposition 2.1} Let $w\in H^0(X,S^2\Omega^1_X)$ be of rank 2 on a $X$  a connected complex surface, then $w$ being closed is equivalent to:

$$K_\Bbb C(w)=0$$

\noindent Moreover, $K_\Bbb C(w)=0$ is equivalent to $w$ is of the
1st kind on $X\setminus \text {Disc}^0_w$.
\endproclaim

\

\

\

\noindent  {\bf 2.2 Characterization of the breakdown locus $B_w$}

\

In this section $w$ is a closed symmetric 2-differential of rank
2. We start by showing that the breakdown locus $B_w$ has no
isolated points and then proceed to show that $B_w$ is an analytic
subvariety of codimension 1 and to characterize geometrically its
 components.

\proclaim {Lemma 2.2} $B_w$ has no isolated points.
\endproclaim

\demo {Proof} It is enough to show that if $w$ is of the 1st kind
in a punctured ball $\Bbb B^*$, then it is of the 1st kind on the
whole ball $\Bbb B$. Let $w\in H^0(\Bbb B,S^2\Omega^1_\Bbb B)$ be
of the 1st kind on the punctured ball $\Bbb B^*$. According to
corollary 1.2,

$$w|_{\Bbb B^*}=\phi_1\phi_2$$

\noindent with $(\phi_1,\phi_2)\in H^0(\Bbb
B^*,\Omega^1_{X,cl}\otimes (\Bbb C_{\rho}\oplus\Bbb C_{\rho}^*))$.
The triviality of the fundamental group of $\Bbb B^*$ implies that
$\Bbb C_{\rho}\simeq \Bbb C_{\rho}^*\simeq \Bbb C$ and hence the
$\phi_i$ can be chosen to be in $H^0(\Bbb B^*,\Omega^1_{X,cl})$.
Again using $\pi_1(\Bbb B^*)=\{e\}$, it follows by integration
that $\phi=df_i$ with $f_i\in \Cal O(\Bbb B^*)$. Hartog's
extension theorem implies that exist $\hat f_i \in \Cal O(\Bbb B)$
extending the $f_i$ and hence $w$ is of the 1st kind on  $\Bbb B$
with $w=d\hat f_1d\hat f_2$.$\blacklozenge$
\enddemo

Proposition 2.1 tell us that that:

\proclaim {Corollary 2.2} $B_w\subset \text {Supp} (\text
{Disc}^0_w)$.
\endproclaim

\

To proceed we need to give a geometric description of the
irreducible components of both discriminant loci. The support of
the discriminant divisor decomposes into:

$$\text {Supp} (\text {Disc}_w)=N_w\cup S_w$$

\noindent where $N_w$ and $S_w$ are the union of all irreducible
components of $\text {Disc}_w$ of respectively odd and even
multiplicities. The locus $N_w$ corresponds to the points where
$w$ fails to split at. For the support of the core discriminant
divisor we have:

$$\text {Supp} (\text {Disc}^0_w)=N_w\cup C_w \cup R_w$$

\noindent It follows from the characterization of $\text
{Disc}^0_w$ given after definition 2.1 that any $x\in \text {Supp}
(\text {Disc}^0_w)$ is in the closure of the locus of points $y\in
X$ such that $\hat w_i(y)$ defines a single line in $T_yX$ (where
$w|_{U_i}=h_i\hat w_i$ with $y\in U_i$, $h_i\in \Cal O(U_i)$ and
$\hat w_i$ is a symmetric 2-differential vanishing at most in
codimension 2). Note that definition 2.1 gives directly that the
divisor $N_w$ is fully contained in $\text {Supp} (\text
{Disc}^0_w)$, since only even multiples of irreducible components
are subtracted from $(\text {Disc}_w)$ to obtain $ (\text
{Disc}^0_w)$.

\smallskip

 The
divisor $C_w$ consists of the union of all irreducible components
of $\text {Supp} (\text {Disc}_w)$ which are leaves simultaneously
of the two foliations defined by the $\{\hat w_i\}_{i\in I}$ (a
2-differential defines two foliations where it splits). We will
call the irreducible components of $C_w$ the common leaves of $w$.
The divisor $R_w$ consists of all the irreducible components of
$\text {Supp} (\text {Disc}^0_w)$ that are not in $N_w$ or $C_w$.
These will be the components for which at their  general point $x$
the two different foliations given by $\{\hat w_i\}_{i\in I}$
define leaves that are tangent at $x$ but that do not coincide).

\

\proclaim {Theorem 2.4} $B_w=N_w\cup C''_w$, where $C''_w$ is a
union of curves contained in $C_w$.
\endproclaim
\demo {Proof} The locus $N_w$ is contained in $B_w$ since the
differential $w$ splits on any $x\notin B_w$. Set $X'=X\setminus
N_w$, $C_w'=C_w\cap X'$ and $R'_{w}=R_w\cap X'$ and get:

$$\text {Supp} (\text {Disc}^0_w)\cap X'=C_{w}'\cup R_{w}'$$

\noindent The desired result then follows if we show that the
breakdown locus $B_{w|_{X'}}$ is an union of irreducible
components of $C_w'$ (with $C''_w$ being the closure of this union).

\smallskip

By construction $X'$  is the open subset of $X$ where $w$ is
locally split. Hence given any $x\in X'$, there exists an open
neighborhood $U_x$ of $x$ where $w|_{U_x}=\mu_1\mu_2$, $\mu_i\in
H^0(U_x,\Omega^1_{X})$. We can shrink $U_x$ so that we can
decompose $\mu_i=h_i\hat\mu_i$ with $h_i\in \Cal O(U_x)$ and the
$\hat \mu_i \in H^0(U_x,\Omega^1_{X})$ are either nowhere
vanishing or vanish only at $x$. Frobenius' theorem (if
$\hat\mu_i(y)\neq 0$, then $\exists U_y$ open neighborhood of $y$
where $\hat\mu_i=f_idu_i$, $f_i,u_i\in \Cal O(U_y)$), then implies
that the set $S\subset X'$ consisting of the points $x$ where $w$
fails to have a neighborhood $U_x$ where $w|_{U_x}=gdz_1dw_1$ with
$g\in \Cal O(U_x)$ and $dz_1$, $dw_1$ nowhere vanishing is
discrete.

\smallskip

Consider the irreducible decomposition

$$C'_w\cup R'_w=\bigcup_{i=I}C'_{w,i} \cup \bigcup_{j=J}R'_{w,j}$$

\noindent where $I,J$ are countable and $C'_{w,i}$ and $R'_{w,j}$
are the irreducible components of $C'_w$ and $R'_w$ respectively.
Below, we will first show that  the irreducible components
$R'_{w,j}$ intersect $B_{w|_{X'}}$ only inside $S$, i.e.
$R'_{w,j}\cap B_{w|_{X'}}\subset S\cup \bigcup_{i=I}C'_{w,i}$.
Second, we will show that the irreducible components $C'_{w,i}$
are such that $C'_{w,i} \subset B_{w|_{X'}}$ or $C'_{w,i}\cap
B_{w|_{X'}}\subset S$. These two results (and corollary 2.3) imply
that $B_{w|_{X'}}= \bigcup_{i=I'}C'_{w,i} \cup S'$, with
$S'\subset S$ and $I'\subset I$.  The result then follows since
$S'\subset \bigcup_{i=I'}C'_{w,i}$. The discreteness of $S$ and
$\bigcup_{i=I'}C'_{w,i}$ being an analytic subvariety of $X'$
implies that if $x\in S'$ is not contained in
$\bigcup_{i=I'}C'_{w,i}$, then $x$ has a neighborhood $U_x$ such
that $U_x\cap B_w=x$, but by lemma 2.2 $B_w$ has no isolated
points.

\

Claim: $R'_{w,j}\cap B_{w|_{X'}}\subset S \cup
\bigcup_{i=I}C'_{w,i}$

\

Before proceeding, note that by the definition of the set $S$ it
follows that any $x\in X'\setminus S$ has a neighborhood $U_x$
with $g,z_1,z_2,w_1\in \Cal O(U_x)$ such that

$$w|_{U_x}=gdz_1dw_1$$

\noindent and $\phi=(z_1,z_2):U_x \to \Delta\times\Delta$, $\Delta$ a disc centered at $0$, is a biholomorphism with $\phi(x)=(0,0)$.

\smallskip

We will show that  any $x\in R'_{w,j} \cap [X'\setminus (S \cup
\bigcup_{i=I}C'_{w,i})]$ cannot lie in $B_w$.

\smallskip

Let $U_x$ be a neighborhood of $x$ as in the previous paragraph.
Consider the leaf $L=\{z_1=0\}$ of $w$ on $U_x$ passing through
$x$. By hypothesis  $x$ is not in a common leaf of  $w$, hence $L$
can not be a common leaf of $w$ which implies that $L\not \subset
\text {Supp} (\text {Disc}^0_w)$. If $L\subset \text {Supp} (\text
{Disc}^0_w)$ then $dz_1\wedge dw_1=0$ on $L$ making $L$ a leaf of
$dw_1$ also, hence a common leaf for $w$. Hence $L\setminus [\text
{Supp} (\text {Disc}^0_w) \cap L)]\neq \emptyset$

\smallskip

Pick $y\in L$ but not in $\text {Supp} (\text {Disc}^0_w)$, then
by proposition 2.1 there is a (connected) neighborhood $U_y$ of
$y$ where $w|_{U_y}=f(z_1)g(w_1)dz_1dw_1$ with $f,h\in \Cal
O(U_y)$ ($f(z_1)$ denotes a function $f(z_1,z_2)$ depending only
on $z_1$). Let $\Delta'$ be a disc centered at 0 such that
$\Delta'\times z_2(y)\subset \phi(U_y)$ and
$W_x=z_1^{-1}(\Delta')(=\phi^{-1}(\Delta'\times \Delta))$. The
function $f$ has a clear holomorphic extension $\hat f\in \Cal
O(U_y\cup W_x)$, with $\hat f|_{W_x}(z_1,z_2)=f(z_1,z_2(y))$.

\smallskip

The same reasoning applied to $h$ will not give an extension of
$h$ to $W_x\cup U_y$, so instead we use the extension of $f$ and
consider the function $\hat h=\frac {g}{\hat f}$. Clearly, $\hat
h|_{U_y}=h$ hence $\hat h$ is a function of $w_1$ alone. The
function $\hat h|_{W_x}$ is holomorphic since the irreducible
components of the polar divisor $(\hat h|_{W_x})_\infty$ if they
exist must be some of the irreducible components of the divisor of
zeros of $\hat f$ which will be a union of curves $\{z_1=c\}$ and
hence intersect non-trivially $U_y$ but this intersection must be
empty since $h|_{U_y}=h$ is holomorphic.

\smallskip

It follows from the previous two paragraphs that the closed
holomorphic decomposition of $w$ at $U_y$,
$w|_{U_y}=f(z_1)h(w_1)dz_1dw_1$, propagates to give the closed
holomorphic decomposition on the neighborhood $W_x$ of $x$,
$w|_{W_x}=\hat f(z_1)\hat h(w_1)dz_1dw_1$, making $x\not \in B_w$.

\

Claim: $C'_{w,i} \subset B_{w|_{X'}}$ or $C'_{w,i}\cap
B_{w|_{X'}}\subset S$.

\

In addition to the properties, described two paragraphs above,
that we can guarantee for an open neighborhood $U_x$ of $x\in
X'\setminus S$,  we can equally guarantee the existence of an open
neighborhood $U_x'\subset U_x$ and $w_2\in \Cal O(U_x')$ such that
$\phi'=(w_1,w_2):U_x' \to \Delta'\times \Delta'$ is a
biholomorphism.

\

Consider the subsets $C_{w,i}^*=C'_{w,i}\cap (X'\setminus S)$ and
$V_i=C_{w,i}^*\cap (X\setminus B_w)$. The set $C_{w,i}^*$ is
connected since by the local parametrization theorem [De12] an
irreducible component of an analytic variety punctured by a
discrete set is connected. The subset $V_i$ consisting of all
points of $C_{w,i}^*$ where $w$ has a local holomorphic
decomposition is clearly open in $C_{w,i}^*$. We proceed to show
that $V_i$ is also closed in $C_{w,i}^*$. Since $C_{w,i}^*$ is
connected, $V_i$ being both open and closed implies the desired
result that the irreducible components $C'_{w,i}$ are such that
$C'_{w,i} \subset B_{w|_{X'}}$ (when $V_i=\emptyset$ and use $B_w$
closed) or $C'_{w,i}\cap B_{w|_{X'}}\subset S$ (when
$V_i=C_{w,i}^*$).

\smallskip

Let $x\in C_{w,i}^*$ be an accumulation point of $V_i$. Pick $y\in
V_i\cap U_x'$, with $U_x'$ as in two paragraphs above. Since
$C_{w,i}^*$ is a common leaf of $w$, $y\in L_x=\{z_1=0\}
=\{w_1=0\}$. Hence $y$ has a neighborhood $U_y$ such that
$\phi(U_y)\supset \Delta''\times z_2(y)$ and $\phi'(U_y)\supset
\Delta''\times w_2(y)$, $\Delta''$ a disc centered at 0, where
$w|_{U_y}=gdz_1dw_1=f(z_1)h(w_1)dz_1dw_1$ with $f,h\in \Cal
O(U_y)$. The functions $f$ and $h$ are clearly extendable to $\hat
f \in \Cal O(z_1^{-1}(\Delta'')\cup U_y)$ and $\hat h\in \Cal
O(w_1^{-1}(\Delta'')\cup U_y)$. By construction
$W_x=z_1^{-1}(\Delta'') \cap w_1^{-1}(\Delta'')$ is a connected
open set  containing $x$ and $y$ and $g|_{W_x\cap U_y}=\hat f\hat
h|_{W_x\cap U_y}$, hence $g|_{W_x}=\hat f\hat h|_{W_x}$  giving a
holomorphic decomposition of $w$ on the neighborhood $W_x$ of $x$,
i.e. $x\in V_i$.
\enddemo

\

\

\noindent  {\bf 2.3 Monodromy at $B_w$}

\

\

Let $w\in H^0(X,S^2\Omega^1_X)$ be closed of rank 2 and
$B_w=\sum_{j\in J}$, $J$ countable, be the irreducible
decomposition of the breakdown locus $B_w$. Let  $\Cal U=\{U_i\}$
be a Leray covering of $X\setminus B_w$ where

$$w|_{U_i}=\phi_{1i}\phi_{2i}$$

\noindent with $\phi_{ki}\in H^0(U_i,\Omega^1_{X,cl})$, $k=1,2$. The abelian rank of a closed
symmetric 2-differential of rank 2 is trivial, it follows then from proposition 1.1 that if $U_i\cap U_j\neq \emptyset$, then

$$\bmatrix \phi_{1i}\\\phi_{2i}\endbmatrix=g_{ij}\bmatrix \phi_{1j}\\\phi_{2j}\endbmatrix$$

\noindent  with $g_{ij}\in G=\{\bmatrix  c&0\\
 0&c^{-1}\endbmatrix,  \bmatrix  0&c\\
c^{-1}&0\endbmatrix|\text { } c\in \Bbb C^*\}$. The collection $\{g_{ij}\}$ gives a 1-cocycle with
values in the group $G$, i.e. $\{g_{ij}\} \in Z^1(\Cal U,G)$. Hence given $x_0\in X\setminus B_w$,
we obtain a representation $\rho:\pi_1(X\setminus B_w,x_0) \to G$.

\

If $w$ is split, then $\text {Im}\rho \subset G'$, with $G'=\{\bmatrix  c&0\\
 0&c^{-1}\endbmatrix,  \forall c\in \Bbb C^*\}\subset G$ (this follows from being able to get a
 consistent ordering of the foliations on all the $U_i$).  Since $G'$ is abelian we get a
 representation, $\rho_w:\pi_1(X\setminus B_w) \to G'$, that is independent of the base point and
 factors through $H_1(X\setminus B_w,\Bbb Z)$, and gives:

 $$\bar \rho_w:H_1(X\setminus B_w,\Bbb Z) \to G'$$

 \noindent Associated with each irreducible component $B_j$, let $\gamma_ij\in H_1(X\setminus B_w,\Bbb Z)$
 be the class of  a simple loop around $B_ij$  (boundary to a disc transversal to $B_j$ centered at a
 general point of $B_j$) which can have either orientation.

\proclaim {Definition 2.2} Let $w\in H^0(X,S^2\Omega^1_X)$ be
split, closed of rank 2  and $B_w=\sum_{j\in J} B_j$, $J$
countable, be the irreducible decomposition of the breakdown locus
$B_w$. To each irreducible component $B_j$ we associate the
monodromy index $M(B_j,w)=\{c,c^{-1}\}$,
if $\bar \rho_w(\gamma_j) =\bmatrix  c&0\\
 0&c^{-1}\endbmatrix$, with $\bar \rho_w$ and $\gamma_j$ as above.
\endproclaim

\

\

\

\noindent  {\bf 2.4 Local form  at $B_w$}

\

\

The goal of this section and the main result of this article is to give the general form of a split closed symmetric 2-differential of rank 2 $w$ at the general point of an irreducible component of the breakdown locus $B_w$.
We will see that  the closed decompositions of $w$ can acquire essential singularities and have non-trivial monodromy at the breakdown   locus $B_w$. We will show that the essential singularities have an algebraic feature, they come from exponential functions with meromorphic functions with poles along $B_w$ as exponents.  Moreover, we will give a bound on the order of the poles of the meromorphic functions appearing as exponents. The bounds come from the order of contact of the two foliations of $w$ along  the irreducible components of $B_w$.

\

We start with some examples of closed symmetric 2-differentials for which $B_w$ is non-empty.

\

\noindent Example: (non-split) Let $z_1$ be a holomorphic  coordinate of $\Bbb C^n$  and $f\in \Cal O(\Bbb C^n)$,  set $w=z_1(dz_1)^2-(df)^2$ . The differential is non split at all points in $\{z_1=0\}$ but it is closed since any point $y \in X\setminus \{z_1=0\}$  has a  neighborhood $U_y$ where $\sqrt {z_1}$ exists and hence $w$ has a holomorphic exact decomposition $w|_{U_y}=d(\frac {2}{3}z_1^{\frac {3}{2}}+f)d(\frac {2}{3}z_1^{\frac {3}{2}}-f)$.

\

If the differential is locally split at $x$, then  a 2nd layer of the failure of $w$ to have a holomorphic closed decomposition at $x$ is due to the  monodromy in the factors of the closed decompositions (not the  monodromy of the foliations) around $B_w$.

\

\noindent Example (monodromy of the closed decompositions): Let $B\subset \Bbb C^2$ be  a sufficiently small open ball about the origin where $1+z_2$ is invertible. Consider $w=(1+z_2)^\alpha dz_1d[z_1(1+z_2)]$.  Recall that the differential $w$ has a holomorphic closed decomposition at a point $x\in B$ if and only if we can decompose $(1+z_2)^\alpha$ as a product of holomorphic functions of $z_1$ and $z_1(1+z_2)$ near $x$.  At points in the complement of $\{z_1=0\}$ we have the decomposition $(1+z_2)^\alpha=z_1^{-\alpha}[z_1(1+z_2)]^\alpha$, but at points in $\{z_1=0\}$ the functions involved are multivalued, hence no holomorphic closed decomposition of $w$ is possible at $x\in \{z_1=0\}$. In fact this monodromy is infinite if $\alpha\not \in \Bbb Q$, meaning that even after finite ramified coverings the symmetric differential would not have an exact decomposition along the pre-image of $\{z_1=0\}$.

\smallskip

If the differential is both locally split at $x$ and no monodromy occurs, then $w$ a 3rd level of failure to have a holomorphic closed decomposition is due to the singularities of the 1-differentials on the decomposition.

\

\noindent Example: (meromorphic singularities) $w=(dz_1)^2+z_1z_2dz_1dz_2=dz_1(dz_1+z_1z_2dz_2)$ has  the common leaf $L=\{z_1=0\}$. The differential $w$ is closed because the 1-form $dz_1+z_1z_2dz_2$ has an integrating factor, $\frac {1}{z_1}$, which is a function of $z_1$. This integrating factor produces the closed meromorphic decomposition $w=d(\frac {z_1^2}{2})(\frac {1}{z_1}dz_1+z_2dz_2)$. Note that since the abelian rank of $w$ is trivial any other closed decomposition of $w$ would differ just by multiplicative constants hence meromorphic singularities would be always present in the closed decompositions of $w$.

\

\noindent Example: (essential singularities) This example shows that even essential singularities can occur, $w= e^\frac {z_2}{1+z_1z_2}dz_1d[z_1(1+z_1z_2)]$. The 1-differentials in the split closed decomposition are unique up to multiplicative constants, as it was shown in proposition 1.1, and the constants will cancel each other so in fact the decomposition is unique and has the form

$$w=e^\frac {z_2}{1+z_1z_2} dz_1d[z_1(1+z_1z_2)]=e^{\frac {1}{z_1}}dz_1e^{-\frac {1}{z_1(1+z_1z_2)}}d[z_1(1+z_1z_2)]$$

\noindent with essential singularities occurring on the closed 1-forms at $\{z_1=0\}$.

\

\proclaim {Lemma 2.5} Let $X$ be a complex 2-manifold, $w\in H^0(X,S^2\Omega^1_X)$ be split of
rank 2 and $L$ be an irreducible component of a common leaf of $w$. Then there is an $m\in \Bbb N$
such that the  general point $x$ of $L$   has a neighborhood $U_x$ with a holomorphic chart $(z_1,z_2)$ where

$$w|_{U_x}=f(z_1,z_2)dz_1d[z_1(1+z_1^mz_2)] $$

\endproclaim

\demo {Proof} The differential $w$ being split implies that every point $x\in X$ has an open
neighborhood $U_x$ such that $w|_{U_x}=\mu_1\mu_2$, $\mu_i\in H^0(X,\Omega^1_X)$. By shrinking $U_x$
we can factor the 1-forms $\mu_i$ in the form $\mu_i=f_i\hat \mu_i$ with $f_i\in \Cal O(U_x)$ and
$\hat \mu_i$ either non-vanishing or vanishing only at $x$. Since by Frobenius theorem a
non-vanishing 1-form in dimension 2 is integrable, it follows that there is a discrete set
$S\subset X$ such that all $x\in X\setminus S$ have a neighborhood $U_x$ where

$$w|_{U_x}=hdvdr \tag 2.3$$

\noindent with $h, v,r \in \Cal O(U_x)$, $v(x)=r(x)=0$, $dv$ and $dr$ nowhere zero  on $U_x$.

\

Let $x$ be a general point of  $L$, using the notation of (2.3) we have $L\cap U_x=\{v=0\}=\{w=0\}$ with:

 $$r=vu$$

 \noindent with $u$ a unit on $U_x$. After shrinking $U_x$ we can assume there is a holomorphic local chart
 on $U_x$, $(v_1,v_2)$ such that $v_1=v$. Consider the series expansion
 $u(v_1,v_2)=\sum_{i=0,j=0}^\infty c_{ij}v_1^iv_2^j$ and let
 $m=\text {min}\{i| \exists j>0 \text { s.t. } c_{ij}\neq 0 \}$ (the Taylor series of  $u$  must involve
 $v_2$ since $dv\wedge dr\not\equiv 0$). Then decompose $u$ as

 $$u(v_1,v_2)=s(v_1)+v_1^m(t(v_2)+v_1g(v_1,v_2))$$

 \noindent where $s(v_1)=\sum_{i=0}^\infty c_{i0}v_1^i$ is a holomorphic function in $v_1$ with $s(0)\neq 0$
 and hence a unit in a neighborhood of $0$. Note $t(0)=0$  and more importantly $t(v_2)$  is not constant.
 Hence $dt(v_2)$ is non vanishing at the general point of $L\cap U_x$. If $dt(v_2)(x)=0$, then change $x$
 to make $dt(v_2)(x)\neq 0$.

 \

 Set $z_1=v_1s(v_1)$ and $z_2=\frac {t(v_2)+v_1g(v_1,v_2)}{s(v_1)^{m+1}}$. By construction $dz_1(x)\neq 0$, $dz_2(x)\neq 0$ and $dz_1\wedge dz_2(x)\neq 0$ and

 $$r=z_1(1+z_1^mz_2)$$

 \noindent giving the desired $w|_{U_x}=f(z_1,z_2)dz_1d[z_1(1+z_1^mz_2)]$ with $f=\frac {h}{s(v_1)+v_1s'(v_1)}$.

 \

 Observing that $m=\text {ord}_{\{v_1=0\}}(\frac {\partial r}{\partial v_2})-1$, it follows that $m$ is independent of the choice of $v$ and $r$ with $dv$ and $dr$ non-vanishing such that $w=hdvdr$ and the choice of holomorphic chart $(v_1,v_2)$ with $v_1=v$. The independence of $m$ on the above choices plus the connectedness of $L$ minus a discrete set of points implies that any other general point of $L$ would give the same $m$ and hence $m$ is naturally associated to the irreducible component $L$.$\blacklozenge$
 \enddemo

\

\proclaim {Definition 2.3} An irreducible component $L$ of a common leaf  of $w\in H^0(X,S^2\Omega_X^1)$ of rank 2 is said to have order of contact $m$, $O(L,w)=m$, if in a neighborhood $U_x$ of the general point $x\in L$ $w$ is of the form as in lemma 2.5, i.e. $w|_{U_x}=f(z_1,z_2)dz_1d[z_1(1+z_1^mz_2)] $.
\endproclaim

\

\proclaim { Theorem 2.6} Let $X$ be a complex 2-manifold, $w\in H^0(X,S^2\Omega^1_X)$ be split,
closed of rank 2 and $L$ an irreducible component of a common leaf of $w$. Then  the   general
point $x$ in $L$  has a neighborhood $U_x$ where $w|_{U_x}$ has a decomposition of the form:

$$w|_{U_x}=z_1^k(1+z_1^mz_2)^\alpha e^{f(z_1)}e^{g(z_1(1+z_1^mz_2))}dz_1d[z_1(1+z_1^mz_2)]$$

\noindent where:

i) $m=O(L,w)$, $k=\text {ord}_L(w)_0$ ($(w)_0$ is the divisorial zero of $w$) and
 $\alpha={\frac {\log c}{2\pi i}}+k$ for some $k\in \Bbb Z$ and $c\in M(L,w)$.
\smallskip

ii)  $f$ and $g$ are meromorphic functions on $\Delta^*$ with poles of order at most $m$ at $0$.

\endproclaim

\noindent {\bf Remark}:  The local form of $w|_{U_x}$ in the
theorem can be rewritten as the following decomposition of
$w|_{U_x}$ as the product of two closed 1-differentials (in a
generalized sense since they might be multi-valued) with
singularities along $L$:

$$w|_{U_x}=(z_1^\beta e^{f(z_1)} dz_1)([z_1(1+z_1^mz_2)]^\alpha e^{g(z_1(1+z_1^mz_2))}d[z_1(1+z_1^mz_2)])$$

\noindent with $\alpha+\beta=\text {ord}_L(w)_0$ and  $\alpha$, $f$ and $g$ as in the theorem.

\

\demo {Proof} According to the lemma 2.5  the general point $x\in
L$ has a neighborhood $U_x$ with a holomorphic coordinate chart
$(z_1,z_2)$ such that $x=(0,0)$, $L\cap U_x=\{z_1=0\}$ and
$w|_{U_x}=v(z_1,z_2)dz_1d[z_1(1+z_1^mz_2)]$ with $v \in \Cal
O(U_x)$.

\smallskip

We claim that if we shrink $U_x$ the divisor of zeros of
$w|_{U_x}$ is $(v)_0=kL$ and hence

$$w|_{U_x}=z_1^k\tilde w$$

\noindent with $\tilde w\in H^0(U_x,S^2\Omega^1_X)$ a nowhere
vanishing closed symmetric differential of the form:

$$\tilde w=\tilde v(z_1,z_2)dz_1d[z_1(1+z_1^mz_2)] \tag 2.4$$

\noindent with $\tilde v(z_1,z_2)\in \Cal O^*(U_x)$.

\smallskip

By shrinking $U_x$ we can make $(v)_0$ a finite union of
irreducible components all passing through $x$. The differential
$w$ being closed implies (theorem 2.4) that all $y\in U_x\setminus
L$ have a neighborhood $U_y$ such that
$v|_{U_y}=f(z_1)g(z_1(1+z_1^mz_2))$. This implies that if an
irreducible component of $(v)_0$ is not $L$, then it must be a
level set of $z_1$ or $z_1(1+z_1^mz_2)$ not passing through $x$, a
contradiction. It follows then that $(v)_0=kL$ for some $k\in \Bbb
N$ and (2.4) holds.

\smallskip

 Note that we
have the equality  $M(L,\tilde w)=M(L,w)$, this can be seen for
example by noting that the factor on the local holomorphic
decompositions of $w$ and $\tilde w$ corresponding to the
foliation $d[z_1(1+z_1^mz_2)]$ does not change (the 1-cocycle with
values in $\Bbb C^*$ corresponding to this foliation remains
unchanged) hence the index remains unchanged.

\

The neighborhood $U_x$  can be chosen to be the bi-disc
$U_x=\Delta_{\epsilon_1}\times \Delta_{\epsilon_2}$,
$\epsilon_i>0$, relative to the coordinate chart $(z_1,z_2)$. On
$U_x$ we have two maps $\pi_1:U_x \to \Bbb C$ given by
$\pi_1(z_1,z_2)=z_1$ and $\pi_2:U_x \to \Bbb C$ given by
$\pi_2(z_1,z_2)=z_1(1+z_1^mz_2)$.

\

Let $\Cal U=\{U_i\}_{i=1,...,k}$, $k\in \Bbb N$, be a Leray  open covering of the punctured disc $\Delta_{\epsilon_1}^*$.
The Leray covering $\{U_i\times \Delta_{\epsilon_2}\}_{i=1,...,k}$ of $U_x\setminus \{z_1=0\}$ is such that one has the
holomorphic closed decompositions on its open sets:

$$\tilde w|_{U_i\times \Delta_{\epsilon_2}}=\breve f_i(z_1)\breve g_i(z_1(1+z_1^mz_2))dz_1d[z_1(1+z_1^mz_2)] \tag 2.5$$

\noindent where $\breve f_i=\pi_1^*f_i$ with $f_i\in \Cal O(U_i)$
and $\breve g_i=\pi_2^*g_i$ with $g_i\in \Cal O(U'_i)$
$(U_i\subset U'_i=\pi_2(U_i\times \Delta_{\epsilon_2}))$. The
existence of such closed decomposition on the whole open sets
$U_i\times \Delta_{\epsilon_2}$ is guaranteed since the open sets
are simply connected and the fibers of both $\pi_1|_{U_i\times
\Delta_{\epsilon_2}}$ and $\pi_2|_{U_i\times \Delta_{\epsilon_2}}$
are connected (assuming $\epsilon_1$ and $\epsilon_2$ are
sufficiently small).

\

Since $w$ is a symmetric differential of degree 2 and rank 2, the
abelian rank of $w$ is trivial which due to proposition 1.1
implies that on the intersections $U_i\cap U_j$:

$$ f_i=c_{ij} f_j \text { }\text { }\text { }\text { }\text { }\text { }\text { }\text { } g_i=c_{ij}^{-1} g_j \tag 2.6$$

\noindent The collection $\{c_{ij}\}$ defines a 1-cocycle in
$Z^1(\Cal U,\Bbb C^*)$ defining a representation of $\rho:
\pi_1(\Delta_{\epsilon_1}^*) \to \Bbb C^*$.  There is a natural
homomorphism $\phi_L:\pi_1(\Delta_{\epsilon_1}^*)\to
H_1(X\setminus B_w,\Bbb Z)$ sending  the class of a simple loop
$\gamma$ around the origin oriented positively to  the class of a
simple loop $\gamma_L$ around the irreducible component $L$ (as in
definition 2.2). By construction, $\rho(\gamma)$ is one of the
diagonal entries of $\bar \rho_w(\gamma_L)$, i.e. $\rho(\gamma)\in
\text {M(L,w)}=\{c,c^{-1}\}$.

\

To simplify notation rescale the coordinates so that
$U_x=\Delta\times \Delta$, $\Delta$ the unit disc centered at $0$
and set the covering $\Cal U$ of $\Delta^*$ to be
$\{U_{-1},U_0,U_1\}$ with $U_i=(0,1)\times
(\frac{(2i-1)}{3}\pi-\epsilon,\frac {(2i+1)}{3}+\epsilon)$,
$\epsilon>0$ sufficiently small, if expressed in polar
coordinates.

\

Consider the universal covering map $e:\Cal H^- \to \Delta^*$,
$z\to e^z$, with $\Cal H^-=\{z\in \Bbb C|\text {Re} z<0\}$, and
the open covering of $\Cal H^-$, $\tilde \Cal U=\{\tilde
U_j\}_{j\in \Bbb Z}$ where the $\tilde U_j=(-\infty,0)\times
(\frac{(2j-1)}{3}\pi-\epsilon,\frac {(2j+1)}{3}+\epsilon)$. Note
that $e:\tilde U_j\to U_{[j]}$, with $[j]\in \{-1,0,1\}$ and
$j\equiv [j]$ mod 3, is a biholomorphism.

\

Let $\{\tilde f_j\}_{j\in \Bbb Z}\in C^0(\tilde \Cal U,\Cal O^*)$
be the 0-cochain defined by $\tilde f_j=f_{[j]}\circ e \in \Cal
O^*(\tilde U_j)$. The co-boundary $\delta\{\tilde f_j\}$ gives a
1-cocycle with values in $\Bbb C^*$,   $\{\tilde c_{jj'}\}\in
Z^1(\tilde \Cal U,\Bbb C^*)$, since $\tilde f_j=\tilde
c_{jj'}\tilde f_{j'}$ on $\tilde U_j\cap \tilde U_{j'}$. The space
$\Cal H^-$ being simply connected implies that $\{\tilde
c_{jj'}\}\in B^1(\tilde \Cal U,\Bbb C^*)$. Hence there is a
collection $\{\tilde c_j\}\in C^0(\tilde \Cal U,\Bbb C^*)$ such
that $\tilde c_{j'}\tilde f_{j}=\tilde c_j\tilde f_j$ on $\tilde
U_{j}\cap \tilde U_{j'} $ giving:

 $$\{\tilde c_j\tilde f_j\}=:F\in \Cal O^*(\Cal H^-) \tag 2.7$$

 \noindent The function $F$,  due to $\tilde c_{jj'}=c_{[j][j']}$ and
 the discussion following (2.6), satisfies the special transformation law

 $$F(z+2\pi i)=\rho(\gamma)F(z)$$

 \noindent Since $e^{(\frac {\log \rho(\gamma)}{2\pi i})z}$ is
 function with the same transformation law as $F$, it follows that:

 $$F=e^{(\frac {\log \rho(\gamma)}{2\pi i})z}\hat f(e^z)$$

 \noindent  with   $\hat f\in \Cal O^*(\Delta^*)$.

 \

 The above implies that if we set $c_i=\tilde c_i$, $i=-1,0,1$,
 then

  $$c_i f_i|_{\hat U_i}=(z^{\frac {\log \rho(\gamma)}{2\pi i}}\hat f)|_{\hat U_i}$$

  \noindent  with $z^{\frac {\log \rho(\gamma)}{2\pi i}}$
 representing the principal branch of the  power function  and $\hat
U_i=(0,1)\times (\frac{(2i-1)}{3}\pi,\frac {(2i+1)}{3}\pi)$ (if
expressed in polar coordinates). The same reasoning can be done
with respect to the collection
 $\{ g_i\}_{i=-1,0,1}$ using the collection $\{c_i^{-1}\}_{i=-1,0,1}\in C^0(\Cal U,\Bbb C^*)$ to
  obtain

  $$c_i^{-1}\ g_i|_{\hat U_i}=(z^{-\frac {\log \rho(\gamma)}{2\pi i}}\hat g)|_{\hat U_i}$$

  \noindent with $\hat g\in \Cal O^*(\Delta^*)$ and $z^{-\frac {\log \rho(\gamma)}{2\pi i}}$
  representing the principal branch of  power function.

 \

 Finally, using the above descriptions of the collections $\{c_i\ f_i\}$ and $\{ c^{-1}_i g_i \}$
 it follows that we can rewrite the local holomorphic closed
 decompositions described in (2.5) by changing $\breve f_i$ and $\breve g_i$ to
respectively $c_i\breve f_i$
 and $c^{-1}_i\breve g_i$ and obtain the global decomposition of
 $\tilde w$ on $\Delta^*\times \Delta$:

 $$\tilde w|_{\Delta^*\times \Delta}=(1+z_1^mz_2)^{-\frac {\log \rho(\gamma)}{2\pi i}}\hat f(z_1)\hat g(z_1(1+z_1z_2))dz_1d[z_1(1+z_1^mz_2)] \tag 2.8$$

\noindent  recall that by construction $\rho(\gamma) \in M(L,w)$.
Note that behind the global decomposition (2.8) there is a closed
decomposition of $\tilde w$ on $U_x$ but it involves the
multi-valued functions and functions with singularities along $L$,
$\tilde w= (z_1^{\frac {\log \rho(\gamma)}{2\pi i}}\hat
f(z_1)dz_1)([z_1(1+z_1^mz_2)]^{-\frac {\log \rho(\gamma)}{2\pi
i}}\hat g(z_1(1+z_1z_2))d[z_1(1+z_1^mz_2)])$.

\

The next goal is to understand the singularities that are possible
for the functions $\hat f,\hat g\in \Cal O^*(\Delta^*)$. To achieve this goal,
we use the fact that the
product of $\hat f(z_1)$ with $\hat g(z_1(1+z_1^mz_2))$ extends to
a holomorphic function on $\Delta\times \Delta$ since it satisfies:

$$\hat f(z_1)\hat g(z_1(1+z_1^mz_2))|_{\Delta^*\times \Delta}=\hat v(z_1,z_2)|_{\Delta^*\times \Delta} \tag 2.9$$

\noindent where $\hat v(z_1,z_2)=\tilde
v(z_1,z_2)(1+z_1^mz_2)^{\frac {\log \rho(\gamma)}{2\pi i}} \in
\Cal O^*(\Delta\times \Delta)$.

\

The functions $\hat f$ and $\hat g$ do not necessarily have a well
defined logarithm on $\Delta^*$, since $f_*,g_*:\pi_1(\Delta^*)
\to \pi_1(\Bbb C^*)$ are not necessarily trivial. However, if we
set  $k_1=f_*(\gamma)\in \pi_1(\Bbb C^*)=\Bbb Z$  and
$k_2=g_*(\gamma)\in \pi_1(\Bbb C^*)=\Bbb Z$ with $\gamma$ a simple
loop around 0 positively oriented, then $\hat f(z)=z^{k_1} \breve
f(z)$ and $\hat g(z)=z^{k_2}\breve g(z)$ are such that the
functions $\breve f,\breve g \in \Cal O^*(\Delta^*)$ have well
defined logarithmic functions,  $f=\log \breve f,g=\log \breve
g\in \Cal O(\Delta^*)$.

\

It follows from (2.9) that $\hat f(z)\hat g(z)=\hat v(z,0)$ and
hence $\hat f(z)\hat g(z)\in \Cal O(\Delta^*)$ extends to a
holomorphic function on $\Delta$ which forces $k_2=-k_1$ ($(\hat
f\hat g)_*:\pi_1(\Delta^*) \to \pi_1(\Bbb C^*)$ is trivial since
it factors through $\hat v(z,0)_*:\pi_1(\Delta) \to \pi_1(\Bbb
C^*)$ and $(\hat f\hat g)_*(\gamma))=f_*(\gamma)+g_*(\gamma)$).
This implies that decomposition (2.8) can be rewritten as:

$$ \tilde w=(1+z_1^mz_2)^{-\frac {\log \rho(\gamma)}{2\pi i}-k_1} e^{f(z_1)}e^{g(z_1(1+z_1z_2))}dz_1d[z_1(1+z_1^mz_2)] \tag 2.10$$

 We are now interested in the singularities of $f,g\in \Cal
O(\Delta^*)$. It follows from (2.9) that:

$$f(z_1)+ g(z_1(1+z_1^mz_2))=\log \hat v(z_1,z_2) \tag 2.11$$

\noindent where $\log \hat v(z_1,z_2)\in \Cal O(U_x)$. To derive
conditions on $f, g \in \Cal O(\Delta^*)$ from (2.11) consider the
Laurent series expansions:

$$f(z_1)=\sum^\infty_{i=-\infty}a_iz_1^i$$

$$ g(z_1(1+z_1^mz_2))=\sum^\infty_{i=-\infty}b_i[z_1(1+z_1^mz_2)]^i$$

\noindent for the sum $f(z_1)+ g(z_1(1+z_1^mz_2))$ to be holomorphic we must have

$$\sum^{-1}_{i=-\infty}a_iz_1^i+\sum^{-1}_{i=-\infty}b_i[z_1(1+z_1^mz_2)]^i=:r(z_1,z_2) \tag 2.12$$

\noindent with $r(z_1,z_2)$ holomorphic. To simplify our notation, we quickly note that for $r(z_1,0)$ to be holomorphic we must have $b_i=-a_i$ $\forall i<0$ from which it follows that

$$r(z_1,z_2)=\sum^{-1}_{i=-\infty}a_iz_1^i[1-(1+z_1^mz_2)^i]$$

\noindent Consider the expansion $1-(1+z_1^mz_2)^i=\sum_{k=1}^\infty c_k^{(i)}z_1^{km}z_2^k$. Of the coefficients $c_k^{(i)}$ we will
only use the fact that they are all non-vanishing and $r(z_1,z_2)=\sum^{-1}_{i=-\infty}\sum_{k=1}^\infty a_ic_k^{(i)}z_1^{i+km}z_2^k$.

\

Reorganizing the terms of the last expansion of $r(z_1,z_2)$, one obtains:

$$r(z_1,z_2)=\sum^{\infty}_{j=-\infty}(\sum_{k=\text {min}\{1,\ulcorner \frac {j}{m}\urcorner\}}^\infty a_{j-km}c_k^{(j-km)}z_2^k)z_1^j$$

\noindent The holomorphicity of $r(z_1,z_2)$ implies that $\forall j\le -1$ the functions

$$s_j(z_2)=\sum_{k=1}^\infty a_{j-km}c_k^{(j-km)}z_2^k$$

\noindent  must vanish, which using the non vanishing of the $c_k^{(i)}$ implies that

$$a_i=0 \text { }\text { }\text { }\text { }\text { }\text { }\text { }\text { }\text { }\text { }\text { }\text { }\text { } \forall i< -m$$

\noindent this jointly with the equality $b_i=-a_i$ $\forall i<0$ give the desired result ii) stating that $f$ and $g$ are meromorphic functions with poles of order at most $m$ at the origin.$\blacklozenge$
\enddemo

\end